\documentclass[12pt]{article}

      \usepackage{latexsym}
         \usepackage[reqno, namelimits, sumlimits]{amsmath} 
         \usepackage{amssymb, amsfonts}
         \usepackage{amsthm} 

\newtheorem{Theorem}{Theorem}
\newtheorem{Lemma}{Lemma}
\newtheorem{Remark}{Remark}

\newcommand{\be}{\begin{equation}}
\newcommand{\ee}{\end{equation}}

\newcommand{\rr}{\mathbb R}
\newcommand{\rt}{{\rr}^3}
\newcommand{\bro}{B_{R_0}}

\newcommand{\rtsmb}{\rt\setminus{\bar B}_{R_0}}
\newcommand{\tu}{\tilde u}
\newcommand{\tp}{\tilde p}
\newcommand{\cc}{C_*}
\renewcommand{\th}{\theta}
\newcommand{\vf}{\varphi}
\newcommand{\ul}{U_{\mathrm{lin}}}
\newcommand{\pl}{P_{\mathrm{lin}}}
\newcommand{\ve}{\varepsilon}
\newcommand{\tps}{\tilde\psi}
\newcommand{\tpsar}{\tps_{A,r_0}}
\newcommand{\tuu}{\tilde U}
\newcommand{\tuuar}{\tuu_{A,r_0}}
\newcommand{\tpp}{\tilde P}
\newcommand{\tppar}{\tpp_{A,r_0}}
\newcommand{\ub}{U^b}
\newcommand{\pb}{P^b}
\newcommand{\tub}{{\tilde U}^b}
\newcommand{\tpb}{{\tilde P}^b}
\newcommand{\tf}{\tilde F}
\renewcommand{\aa}{\alpha}
\newcommand{\xa}{X_\aa}
\newcommand{\norm}[1]{[\hskip-2.6pt|#1|\hskip-2.6pt ]_\aa}
\newcommand{\rtsmo}{\rt\setminus\{0\}}

\renewcommand{\div}{\operatorname{div}}

\newenvironment{myproof}{\begin{proof}[{\rm\bf Proof.}]}{\end{proof}}
\newenvironment{myproofb}{\begin{proof}[{\rm\bf Proof of Theorem 1.}]}{\end{proof}}

\title{
On the large-distance asymptotics of steady state solutions of the Navier-Stokes equations in
3D exterior domains}

 \author{A. Korolev\footnote{University of Minnesota}
\and V. \v Sver\'ak\footnote{University of Minnesota. Supported in part by NSF
  Grant DMS-0457061}}

\date{November 2, 2007}
     
\begin{document}

\numberwithin{equation}{section}

\maketitle


\begin{abstract} We identify the leading term describing the behavior at large distances of the steady state 
solutions of the Navier-Stokes equations in 3D exterior domains with vanishing velocity at the spatial infinity.
\end{abstract}

\section{Introduction}\label{intro}

We consider the 3D steady-state Navier-Stokes equations in exterior domains and study the behavior of the
solutions ``near infinity".  The equations are
\be
\label{nse}
\left.\begin{array}{rcl}
-\Delta u +u\nabla u +\nabla p & = & 0 \,\,,\\
\div u & = & 0
\end{array}\right\} \qquad \mbox{in $\rtsmb$,}
\ee
where $B_{R_0}$ denotes the ball of radius $R_0$ centered at the origin. 
Our main assumption 
about the solutions will be the decay condition 
\be
\label{condition1}
|u(x)|\le {\frac{\cc}{R_0+|x|}}\qquad \mbox{in $\rtsmb$}\,
\ee
for sufficiently small $\cc$.
The specific boundary conditions at $\partial B_{R_0}$ will play no role in our results. 

Naively one might think that the behavior near infinity of the above solutions should be given by the linearized
equation. An immediate well-known objection\footnote{The objection was probably raised already 
by the classics in the 19th century.}  to that is that we expect decay $|\nabla^k u|=O(|x|^{-k-1})$  and $|\nabla^k p|=O(|x|^{-k-2})$
as $x\to\infty$, and therefore the non-linear term in the equation should have the same order of magnitude as the linear 
terms, making the accuracy of the linearization questionable. This heuristics is made rigorous in~\cite{GaldiDeuring}, where it is proved that the leading order term describing the behavior of the solutions
cannot be given by the linearized equation. 

In this paper we identify explicitly the leading order behavior of the above solutions near infinity. We show that
it is given by the explicit solutions calculated by L. D. Landau  in 1943, see~\cite{Landau, LL}. These calculations
were revisited and certain extensions were obtained in~\cite{Cannone, TX}. The Landau solutions
were recently characterized in \cite{Sverak} as the only solutions of the steady Navier-Stokes equation in
$\rt\setminus\{0\}$ which are smooth and $(-1)$-homogeneous in $\rt\setminus\{0\}$. The Landau solutions in $\rt\setminus\{0\}$
can be parametrized by vectors $b\in\rt$ in the following way: For each $b\in\rt$ there exists a unique $(-1)$-homogeneous
solution $U^b$ of the steady Navier-Stokes equations together with an associated pressure $P^b$ which is $(-2)$-homogeneous, such that
$U^b, P^b$ are smooth in $\rtsmo$, $U^b$ is weakly div-free across the origin
and satisfies
\be
\label{landaudef}
-\Delta U^b + \div (U^b\otimes U^b) + \nabla P^b = b \delta(x)\,\,
\ee    
in the sense of distributions.
Here $\delta(x)$ denotes the Dirac function and we use the standard notation $u\otimes v$ for the tensor field $u_i v_j$ defined by the tensor product of the vector fields 
$u=(u_1,\dots,u_n)$ and $v=(v_1,\dots,v_n)$. We also use the standard notation $\div T$ for the vector
field ${\frac{\partial}{\partial x_j}} T_{ij}$.  (The uniqueness of $U^b, P^b$ is much easier to prove if we add the requirement that
$U^b$ be axi-symmetric with respect to the axis passing through the origin in the direction of the vector $b$, but 
as was shown in \cite{Sverak}, this additional symmetry assumption is not necessary.) 
As noticed by Landau, the solutions can be calculated  explicitly in terms of 
elementary functions, see formulae are~\eqref{landausol},\eqref{lforce} and~\eqref{b}.  (It was observed in \cite{Sverak}
that the Landau solutions are in a natural one-to-one correspondence with the group of conformal transformations
of the two-dimensional sphere, and Landau's formulae can be also derived from this observation by using some
standard geometry.)

If $u,p$ is a solution of~\eqref{nse}, we will denote by 
\be
\label{mflux}
T_{ij}=T_{ij}(u,p) = p\delta_{ij}+u_iu_j-({\frac{\partial u_i}{\partial x_j}}+{\frac{\partial u_j}{\partial x_i}})
\ee
the momentum flux density tensor in the fluid. Our main result is the following:

\begin{Theorem}
\label{mainth}
For each $\aa\in(1,2)$ there exists $\ve=\ve(\alpha)>0$ such that the following statement holds true:
 Let $u,p$ be a solution of\eqref{nse} in $\rt\setminus{\bar B}_{R_0}$ satisfying condition~\eqref{condition1} with $\cc\le\ve$. 
 Let $b=(b_1,b_2,b_3)$ be defined by
 \be
 \label{formulab}
 b_i=\int_{\partial B_{R_1}} T_{ij}(u,p)n_j(x)
 \ee
 for some $R_1>R_0$. (Note that the integral is independent of $R_1$.) Let $U^b$ be the Landau solution corresponding to 
 the vector $b$.
 Then
 \be
 \label{asymptotics} 
u(x)  =  U^b(x)+O(|x|^{-\aa}) 
 \quad\qquad\qquad \mbox{as $|x|\to\infty$}\,\,
 \ee
and, for a suitable constant $p_0$,
\be
\label{asymptoticsp}
p(x)-p_0=P^b(x)+O(|x|^{-\aa-1})\qquad \mbox{as $|x|\to\infty$.}
\ee
\end{Theorem}

\begin{Remark}
{\rm
Standard estimates for the linear Stokes system (such as estimate~\eqref{stokesest} in the next section), together with the scaling 
symmetry $u(x)\to\lambda u(\lambda x)$ of 
Navier-Stokes can be used to show that any solution of~\eqref{nse} satisfying~\eqref{asymptotics} will also satisfy
\be
\label{asymptoticshd}
\nabla^k u(x)=\nabla^k U^b(x)+O(|x|^{-k-\aa})\quad\mbox{as $|x|\to\infty$, for $k=1,2,\dots$,}
\ee
and
\be
\label{asymptoticsphd}
\nabla^k p=\nabla^k P^b + O(|x|^{-k-1-\aa})\qquad\mbox{as $|x|\to\infty$, for $k=1,2,\dots\,\,.$}
\ee
See for example~\cite{SverakTsai} for an argument of this type.
}
\end{Remark}

The existence of expansions similar to~\eqref{asymptotics} with $U^b$ replaced by a less specific term, namely a $(-1)-$ homogeneous
function, was is studied in \cite{NazarovPileckas}. The main result of that paper is, roughly speaking, that under smallness
conditions similar to~\eqref{condition1}, the solutions of~\eqref{nse} are ``asymptotically $(-1)$-homogeneous".
As is shown in \cite{Sverak}, the leading term of the asymptotical expansion at $\infty$ of any solution of~\eqref{nse} which
is  asymptotically $(-1)$-homogeneous must be given by a Landau solution. (This result remains true even for large data, since
the proof is not based on perturbative arguments.) Therefore the results in~\cite{NazarovPileckas} together with the results 
in~\cite{Sverak} imply a version of Theorem~\ref{mainth}. Our proof in this paper is much simpler than the proof one could get by 
combining~\cite{NazarovPileckas} and~\cite{Sverak}.
Also, if one tried to evaluate explicitly the values of the constants in the smallness conditions,
the constants coming from the proof here would probably be more favorable.

 The proof of the Theorem~\ref{mainth} is given in the following sections. The main idea of the proof is as follows.
 Assume for simplicity that $u$ satisfies the ``no outflow to infinity condition"
 \be
\label{nooutflow1}
\int_{\partial B_{R_1}} u(x)\cdot n(x) = 0 
\ee
for some $R_1>R_0$, where $n(x)$ denotes the unit normal to $\partial B_{R_1}$. (Since $u$ is div-free, the last integral
does not depend on $R_1$.)
 We extend the fields $u,p$ to fields defined in all $\rt$ and satisfying the inhomogeneous equation
 \be
\label{nsef}
\left.\begin{array}{rcl}
-\Delta u +u\nabla u +\nabla p & = & f \,\,,\\
\div u & = & 0
\end{array}\right\} \qquad\mbox{in $\rt$.}
\ee
The extension needs to be done in a way which enables one to control smooth norms of the extended
function by the corresponding norms of the original function. 
We have $b=\int_{\rt} f$. We then search for solutions of~\eqref{nsef} in the form
$u={\tilde U}^b+v$, where ${\tilde U}^b$ is a suitable regularization of the Landau
solution $U^b$. The equation for $v$ is  solved (for small data) by a standard 
perturbation analysis in the space of continuous functions with decay $O(|x|^{-\aa})$
as $|x|\to\infty$. For this argument one needs both an existence result (for $v$) and a 
uniqueness result (for $u$),
and therefore it looks unlikely that our method could be used in a large data situation.

The general situation when the outflow $\int_{\partial B_{R_1}}u\cdot n$ does not vanish can be handled 
by a standard method of writing $u=a+v$ where $v$ has no outflow, and
$a$ is a suitable multiple of the canonical outflow field $\frac{x\,\,}{|x|^3}$.
See for example~\cite{Finn}, Section 2.2, \cite{Galdi}, Chapter IX, or \cite{NazarovPileckas}, Remark 3.2. 
Roughly speaking, the part of the flow
which produces a non-zero outflow has decay $O(|x|^{-2})$, and therefore it does not influence the
main term in~\eqref{asymptotics} at large distances. See Section~\ref{perturbation} for details.

One can see easily by looking at the Landau solutions that the best possible decay rate
$\aa$ for which the result might still be true is $\aa=2$. We conjecture
that the result indeed remains true for $\aa=2$. However, as the example in Remark~\ref{remark1}
shows, one would probably need to go beyond the elementary perturbation theory used in this paper 
to prove that.

The problem of steady-state solutions of the Navier-Stokes equations in exterior 3D domains has a long
history going back to Leray's paper \cite{Leray}. Leray proved the existence of solutions 
with finite energy. Such solutions are easily seen to be smooth since the steady state equation
is subcritical with respect to the energy estimate. However, the precise behavior of these solutions
as $|x|\to\infty$ is a more subtle problem. This problem shares some features with the (super-critical) 
regularity problem for the time-dependent equation, since there seems to be some vague duality between
regularity (or short-distance behavior) of super-critical problems and asymptotics at large times/distances
(or long distance behavior) of sub-critical problems. In particular, in both cases it seems to be important to obtain
some local control of the energy flux which is stronger than what one can immediately get from the
known conservation laws.  In this paper we will not address these difficult issues, which arise for 
large data, and we will only treat the small data situation, which can be handled by a simple perturbation
theory, and is independent of the energy methods. For the steady state exterior problem this approach was pioneered
by Finn, see e.\ g.\ \cite{Finn, Galdi}. We note that the 3D exterior problem with non-zero velocity at $\infty$ has been
more or less fully solved, even for large data, see~\cite{Babenko, Galdi}, since the non-zero velocity at infinity
sufficiently regularizes the flow. In the 2D situation many problems remain open
even in the case of non-zero velocity at infinity, see~\cite{Amick, Galdi}.

How reasonable are our assumptions? We note that Finn proved in \cite{Finn} the existence of solutions satisfying
our assumptions for quite general boundary-value problems in exterior domains under smallness assumptions on the data.
See also \cite{Galdi} for extensions of these results (still under smallness assumptions).
It is quite conceivable that Theorem~\ref{mainth} remains true 
even without assuming that the constant $\cc$ in~\eqref{condition1} is small. The proof of such a result would
however require to go beyond the perturbation theory and the standard energy methods.

The following example, taken from \cite{SverakPlechac}, shows that the question of relaxing the decay condition~\eqref{condition1}
to a slower decay might be quite subtle. Consider the equation
\be
\label{mNSE}
-\Delta u + (1-a)u\nabla u+{\frac{a}2}\nabla |u|^2 + {\frac12}u\div u = 0\qquad\mbox{in $\rt$}\,\,,
\ee 
for vector fields $u$ in $\rt$. The number $a\in(0,1)$ is a parameter. (For \hbox{$a=\frac12$} the non-linear term in~\eqref
{mNSE} can be written as $\div Q(u,u)$
for a suitable quadratic expression $Q$.) The equation has the same energy estimate as the Navier-Stokes equations.
It turns our that~\eqref{mNSE} has a nontrivial global smooth solution $\bar u$
satisfying $|\bar u(x)|\sim |x|^{-2/3}$ as $x\to\infty$. Since it appears that the various perturbation and energy methods 
used for the steady Navier-Stokes should also work for~\eqref{mNSE}, at least in the case \hbox{$a=\frac12$}, the properties of
$\bar u$ indicate some limitations to these methods. On the other hand, we should remark that steady Navier-Stokes does
have some special properties which are probably not shared by~\eqref{mNSE}. For example, for steady Navier-Stokes
the quantity ${\frac12}|u|^2+p$ satisfies a maximum principle, see e.\ g.\ \cite{Amick}.

The paper is organized as follows: In Section~\ref{preliminaries} we explain the material necessary
for extending the solutions to $\rt$ in a controlled manner. In Section~\ref{landau} we recall the 
necessary facts about Landau's solutions. Finally, in Section~\ref{perturbation} we explain
the perturbation argument.

\vskip.3truecm
\noindent
\emph{Acknowledgement.} The authors are indebted to Professor G. P. Galdi for his valuable comments on a preliminary 
version of the manuscript, which lead to significant improvements.

\section{Preliminaries}\label{preliminaries}

We consider the solutions of the steady Navier-Stokes equation
\be
\label{nse1}
\begin{array}{rcl}
-\Delta u +u\nabla u +\nabla p & = & 0 \,\,,\\
\div u & = & 0
\end{array}
\ee
which are defined ``in the neighborhood of infinity", i.\ e.\ in the region 
$\rt\setminus{\bar B}_{R_0}$, where $\bro$ denotes the ball of radius $R_0$ 
centered at the origin.  In this section we will be interested in the solutions which satisfy
\be
\label{decay}
|u(x)|\le \frac{C_*}{R_0+|x|} \qquad \mbox{in $\rt\setminus{\bar B}_{R_0}\,\,,$}
\ee
and the ``no outflow to infinity" condition 
\be
\label{nooutflow}
\int_{\partial B_{R_1}} u\cdot n = 0
\ee
for some $R_1>R_0$. (We note that the integral in~\eqref{nooutflow} is independent
of $R_1$, since $\div u=0$.) 

The constant $C_*$ above will play a special role and we distinguish
it from the ``generic constants" which will be denoted by $c$. If $c$ depends on a parameter
$X$ and we want to emphasize this dependence, we will write $c(X)$ instead of $c$.
The value of $c$ can change from line to line.

  By a solution of~\eqref{nse1} we mean a smooth function vector field
$u$ in $\rt\setminus{\bar B}_{R_0}$ which satisfies~\eqref{nse1} for a suitable $p$. (The pressure
$p$ will be considered only as a ``secondary" variable: Instead of saying
``the solution $(u,p)"$, we can just say ``the solution $u$", with the understanding
that~\eqref{nse1} is satisfied for a suitable $p$.)
Various other notions of solutions are used in the literature (e.\ g.\ weak solutions), but under the 
assumption~\eqref{decay} they all coincide are are equivalent to the one defined above.

One reason that the above way of thinking of $p$ only as an auxiliary variable works quite well is that
the linear steady-state Stokes system 
\be
\label{stokes}
\begin{array}{rcl}
-\Delta u +\nabla p & = & \div f \,\,,\\
\div u & = & 0
\end{array}
\ee
satisfies local elliptic estimates of the form
\begin{eqnarray}
\label{stokesest}
\lefteqn{\hskip-1.8truecm ||\nabla u||_{X(B_{x_0,R})}+
||p-(p)_{B_{x_0,R}}||_{X(B_{x_0,R})}\le} \nonumber \\ & & \hskip1truecm c(R, X)||f||_{X(B_{x_0,2R})}+\tilde c(R,X)||u||_{L^1(B_{x_0,2R})}\,\,,
\end{eqnarray}
where $B_{x_0,R}$ denotes the ball of radius $R$ centered at $x_0$, $\,\,(p)_{B_{x_0,R}}$ is the average of $p$ over the ball $B_{x_0,R}$ and $X$ can be any space in which classical elliptic estimates work, such as an $L^p$-space with
$p\in(1,\infty)$ or a H\"older space. The main point of estimate~\eqref{stokesest} is that there is no
$p$ on the right-hand side. See for example \cite{SverakTsai} for details.

The linear estimate~\eqref{stokesest} combined with the standard bootstrapping and scaling arguments 
(using the scaling symmetry $u(x)\to\lambda u(\lambda x)$) imply that solutions of~\eqref{nse1} satisfying 
estimate~\eqref{decay} with $C_*\le M$ also satisfy
\be
\label{hdecay}
|\nabla^k u(x)|\le c(k,M)\frac{C_*}{(R_{0}+|x|)^{k+1}}\qquad \mbox{in $\rt\setminus B_{2R_0}$ for $k=1,2,\dots$}\,\,.
\ee
\smallskip
We now relate the solutions of~\eqref{nse1} in $\rt\setminus{\bar B}_{R_0}$ to the solutions of 
the equation in $\rt$ with non-trivial right-hand side:
\be
\label{nse2}
\left.
\begin{array}{rcl}
-\Delta u +u\nabla u +\nabla p & = & f \,\,,\\
\div u & = & 0
\end{array}
\right\}
\qquad \mbox{in $\rt$.}
\ee

Let $u$ be a solution of~\eqref{nse1} in $\rtsmb$ satisfying~\eqref{decay} with $C_*\le M$ and let $p$ be the
associated pressure, defined up to a constant. Using~\eqref{hdecay} we see that we can in fact choose a ``normalized" $p$ so 
that, for $C_*\le M$, we have
\be
\label{pressureest}
|\nabla^k p|\le c(k,M){\frac{C_*}{(R_0+|x|)^{k+2}}}\qquad \mbox{in $\rt\setminus B_{2R_0}$ for $k=0,1,2,\dots\,\,$}.
\ee

We can now extend $u,p$ from $\rt\setminus B_{3R_0}$ to $\tilde u, \tilde p$ defined in $\rt$ such that
$\div \tu = 0$ in $\rt$ and
\be
\label{extestu}
|\nabla^k \tu(x)|\le c(k,M){\frac{C_*}{(R_0+|x|)^{k+1}}} \qquad \mbox{in $\rt$ for $k=0,1,2,\dots\,\,,$}
\ee
together with
\be
\label{extestp}
|\nabla^k \tp(x)|\le c(k,M){\frac{C_*}{(R_0+|x|)^{k+2}}} \qquad\mbox{in $\rt$ for $k=0,1,2,\dots$}\,\,.
\ee
The construction of the extension $p\to\tp$ is  standard.
To be able to construct the extension $u\to\tu$, we of course need condition~\eqref{nooutflow}.
With~\eqref{nooutflow} satisfied, the existence of a smooth div-free extension $\tu$ (not necessarily 
satisfying~\eqref{extestu}) is also classical. The construction of a div-free
extension satisfying~\eqref{extestu} can be carried out in many ways. 
One can proceed for example as follows: Let $\eta\colon[0,\infty)\to[0,1]$ be a smooth function
such that $\eta(r)=0$ for $r\le2$ and $\eta(r)=1$ for $r\ge\frac52$, and let $\eta_{R_0}(r)=\eta(\frac{r}{R_0})$.
Now set $\tu=\eta_{R_0}u+v$, where $v$ is a suitable solution of the equation $\div v=-u\nabla\eta_{R_0}=g$ which
is compactly supported in $B_{3R_0}$. The equation $\div v = g$ has of course many compactly supported solutions,
but it is possible to construct a solution operator $S\colon g\to v=Sg$ which has the required 
regularity properties. (Such an operator is sometimes called a Bogovskii operator.) See for example
\cite{Bogovskii} or \cite{Galdi}, Chapter III.3 for details.

We have
\be
\label{newnse}
\left.
\begin{array}{rcl}
-\Delta \tu +\tu\nabla \tu +\nabla \tp & = & f \,\,,\\
\div \tu & = & 0
\end{array}
\right\}
\qquad \mbox{in $\rt$\,,}
\ee
where the right-hand side $f$ is supported in ${\bar B}_{3R_0}$ and satisfies
\be
\label{fest}
|\nabla^k f(x)|\le c(k,M){\frac{C_*}{(R_0+|x|)^{k+3}}}\qquad \mbox{in $\rt$ for $k=0,1,2,\dots\,\,.$}
\ee
Dropping the tildes and changing $R_0$, if necessary, we see that the study of solutions of Navier-Stokes
defined in the neighborhood of infinity and satisfying the ``no outflow" condition~\eqref{nooutflow}
and the growth condition~\eqref{decay}
can be reduced to the study of the solutions $u$ of
the inhomogeneous equation~\eqref{nse2} with $f$ supported in $B_{R_0}$ and satisfying~\eqref{fest}, and
$u$ satisfying~\eqref{decay} globally, with $\cc$ replaced by $c\cc$.

\section{The Landau Solutions}\label{landau}
The Landau solutions are smooth $(-1)$-homogeneous
 solutions of the steady-state
Navier-Stokes equations defined in $\rt\setminus\{0\}$. Under the additional assumption
of axial symmetry, these were first calculated by L.D.Landau in 1943, see \cite{Landau, LL}.
In \cite{Sverak} it was proved that we do not get any new solutions if the assumption
of axial symmetry is dropped. To write down the explicit formulae, we will use the standard polar 
coordinates $r,\theta,\varphi$ defined by
\begin{eqnarray*}
x_1 & = & r\sin\theta\cos\varphi\,, \\
x_2 & = & r\sin\th\sin\vf \,,\\
x_3 & = & r\cos\th\,.
\end{eqnarray*}

The explicit formulae in polar coordinates for the Landau solution $U$ and the corresponding pressure $P$ are as follows:
\be
\label{landausol}
\begin{array}{rcl}
U_r & = & {\frac{2} {r}}[{\frac{A^2-1}{(A-\cos\th)^2}}-1]\,\,,\\
U_\th & = & -{\frac{2\sin\th}{r(A-\cos\th)}}\,\,,\\
U_\vf & = & 0\,\,,\\
P & =& -{\frac{4(A\cos\th - 1)}{r^2(A-\cos\th)^2}}\,\,.
\end{array}
\ee
In the above formulae, $A$ is a parameter satisfying $A>1$. The velocity field $U$ can
also be expressed in terms of the stream function
\be
\label{streamf}
\psi={\frac{2r\sin^2\th}{A-\cos\th}}
\ee 
as
\be
\label{velstream}
\begin{array}{rcl}
U_r & = & {\frac{1}{r\sin\th}}{\frac{\partial \psi}{r\partial\th}}\,\,,\\
U_\th & = & -{\frac{1}{r\sin\th}}{\frac{\partial\psi}{\partial r}}\,\,,\\
U_\vf & = & 0\,\,.\\
\end{array}
\ee
The integral curves of the velocity field $U$ are given the the equations
$\psi={\mathrm{const.\,}}$ and $\vf= {\mathrm{const}}$.

Clearly $U,\,\,U\otimes U$ and $P$ are locally integrable, and a direct calculation (see e.\ g.\ \cite{LL})
gives
\be
\label{lforce}
-\Delta U + \div (U\otimes U) + \nabla P = \beta(A)e_3 \delta(x)\,\,,
\ee
where $e_3$ is the unit vector in the positive $x_3$-direction and 
\be
\label{b}
\beta(A)=16\pi\left(A+{\frac{1}{2}}A^2\log{\frac{A-1}{A+1}}+{\frac{4A}{3(A^2-1)}}\,\right)\,\,.
\ee
It is not hard to check that the function $\beta(A)$ is monotonically decreasing in $(1,\infty)$
and maps this interval onto $(0,\infty)$. In particular, $\beta$ has an inverse function
$\gamma\colon(0,\infty)\to(1,\infty)$.

It is instructive to compare the formula~\eqref{streamf} with the corresponding formula
for the linear Stokes system. Namely, the solution of
\be
\label{stokesfund}
\left.
\begin{array}{rcl}
-\Delta \ul + \nabla \pl & = & e_3\delta(x)\,\,, \\
\div \ul & = & 0
\end{array}
\right\}
\qquad \mbox{in $\rt$\,,}
\ee
satisfying $\ul(x)\to 0$ as $x\to\infty$ is given by the stream function
\be
\label{linstream}
\psi_{\mathrm{lin}}={\frac{1}{8\pi}}r\sin^2\th\,\,.
\ee
We can get another useful comparison if we express the solution of the problem
\be
\label{reynolds}
\left.
\begin{array}{rcl}
-\Delta u_\ve + \ve\div(u_\ve\otimes u_\ve) + \nabla p_\ve & = & e_3\delta(x) \,\,,\\
\div u_\ve & = & 0
\end{array}
\right\}
\qquad \mbox{in $\rt$\,,}
\ee
with the condition $u_\ve(x)\to 0$ as $x\to\infty$ in terms of the Landau solutions.
The formula (for $\ve>0$) is
\be
\label{reynoldssol}
u_\ve={\frac1\ve}U|_{A=\gamma(\ve)}
\ee
and $u_\ve$ is given by the stream function
\be
\label{streamr}
\psi_\ve={\frac{r\sin^2\th}{\ve\gamma(\ve)-\ve\cos\th}}\,\,.
\ee

We will now regularize the Landau solutions near the origin in the following way.
Let $r_0>0$. Consider a smooth function $\rho\colon[0,\infty)\to[0,\infty)$ such that
$\rho(r)=0$ for $r\le r_0$, $\rho(r)=r$ for $r\ge 2r_0$ and the $k-$th derivative
$\rho^{(k)}(r)$ is bounded by $c(k)r^{1-k}$, and define
\be
\label{regstream}
\tps=\tpsar={\frac{2\rho(r)\sin^2\th}{A-\cos\th}}\,\,.
\ee
With the help of $\tps$ we now define the regularized velocity field $\tuu=\tuuar$ by 
the formulae~\eqref{velstream}, with $\psi$ replaced by $\tps$. We also define
\be
\label{regpressure}
\tpp=\tppar  = -{\frac{4\rho(r)(A\cos\th - 1)}{r^3(A-\cos\th)^2}}\,\,.
\ee
It is easy to check that for $A\ge A_0>1$ we have
\be
\label{uregest}
|\nabla^k \tuu|\le{\frac{c(k,A_0)}{A(r_0+|x|)^{k+1}}}\qquad\mbox{in $\rt$ for $k=0,1,\dots$}
\ee
and
\be
\label{pregest}
|\nabla^k\tpp|\le{\frac{c(k,A_0)}{A(r_0+|x|)^{k+2}}}\qquad\mbox{in $\rt$ for $k=0,1,\dots$}\,\,.
\ee

So far we have mostly considered the Landau solution which are axi-symmetric with respect to the
$x_3$ axis. However, it is clear from the above that for each non-zero vector $b\in\rt$ there exist
a unique Landau solution $U^b$ and the associated pressure $P^b$ 
which are axi-symmetric with respect to the axis $\rr\cdot b$ and satisfy
\be
\label{genlandau}
\begin{array}{rcl}
-\Delta\ub+\div(\ub\otimes\ub)+\nabla\pb & = & b\,\delta(x)\,\,,\\
\div \ub & = & 0\,\,.
\end{array}
\ee
We also set $U^0=0$.
For each $\ub,\pb$ the above construction of the regularized solutions gives the regularized fields 
$\tub=\tub_{r_0}$ and $\tpb=\tpb_{r_0}$ which, for $|b|\le M$ will satisfy the estimates
\be
\label{guregest}
|\nabla^k \tub|\le c(k,M){\frac{|b|}{(r_0+|x|)^{k+1}}}\qquad\mbox{in $\rt$ for $k=0,1,\dots$}
\ee
and
\be
\label{gpregest}
|\nabla^k\tpb|\le c(k,M){\frac{|b|}{(r_0+|x|)^{k+2}}}\qquad\mbox{in $\rt$ for $k=0,1,\dots$}\,\,.
\ee

\section{Perturbation Analysis}\label{perturbation}
Let $f$ be a sufficiently regular compactly supported vector field in $\rt$. Let $b=\int_{\rt} f$. Let $r_0=1$
and let $\tuu=\tub=\tub_{r_0}$ and $\tpp=\tpb=\tpb_{r_0}$ be the regularizations of the Landau
solutions $\ub, \pb$ corresponding to the vector $b$ constructed in the previous section. 
We will seek solutions of the steady Navier-Stokes equation
\be
\label{nseagain}
\left.\begin{array}{rcl}
-\Delta u +u\nabla u +\nabla p & = & f\,\,, \\
\div u & = & 0
\end{array}\right\}\qquad \mbox{in $\rt$}
\ee
in the form $u=\tuu+v$. We set
\be
\label{F}
\tf = -\Delta \tuu + \tuu\nabla \tuu+\nabla\tpp\,\, .
\ee
It is easy to check that $\int_{\rt}\tf=b=\int_{\rt} f$.
The equation for $v$ becomes
\be
\label{maineq}
\left.\begin{array}{rcl}
-\Delta v +\tuu\nabla v +v\nabla \tuu+ v\nabla v +\nabla q & = & f-\tf \,\,,\\
\div v & = & 0
\end{array}\right\}\qquad \mbox{in $\rt$\,\,.}
\ee
Let us choose a fixed $\aa\in(1,2)$. We will prove that under some smallness assumptions  
equation~\eqref{maineq} has a unique solution $v$ with decay $O(|x|^{-\aa})$ as $|x|\to\infty$.
An important point is that, by our construction, $\int_{\rt}(f-\tf)=0$.
Using the scaling symmetry, we see that we can assume $r_0=1$ without loss of generality.

Let $G=G_{ij}$ be the Green tensor of the linear Stokes operator. We note that the vector field $G_{i3}$
is given by the stream function~\eqref{linstream}. Another explicit formula for $G$ is
\be
\label{green}
G_{ij}(x)={\frac{1}{8\pi}}(-\delta_{ij}\Delta+{\frac{\partial^2}{\partial x_i\partial x_j}})|x| 
\ee
For our purposes here we will only need the following obvious estimate
\be
\label{greenest}
|\nabla G(x)|\le{\frac c{\,\,|x|^2}}\,\,.
\ee

The required solutions of~\eqref{maineq} will be found for small data by a standard perturbation argument. 
Let $\xa$ be the space of all continuous div-free vector fields $u$ in $\rt$ satisfying $u(x)=O(|x|^{-\aa})$
as $x\to\infty$. A natural norm in $\xa$ is given for example by
\be
\label{norm}
\norm{u}=\sup_{x}\,\,(1+|x|)^\aa\,|u(x)|\,\,.
\ee

Our perturbation analysis is based on the following elementary estimates:
\begin{Lemma}
\label{elementaryest}

Using the notation above, let $b=\int_{\rt} f$ be the vector used in the construction of $\tuu$.
 Then for $|b|\le M$ we have 
 \be
\label{estimate1}
\norm{G*\div(\tuu \otimes v + v\otimes\tuu)}\le c(\aa,M)|b|\norm{v}\,\, ,
\ee
 and
\be
\label{estimate2}
\norm{G*\div(v\otimes w)}\le c(\aa) \norm{v}\norm{w}\,\,.
\ee 
\end{Lemma}

\begin{myproof} The proof these estimates is standard. We move the derivatives to $G$, use the definition of the norm,
after which the only remaining task is to estimate the integral
\be
\label{integral}
I(x)=\int_{\rt}{\frac{dy}{|x-y|^2(1+|y|)^{\aa+\beta}}}\,\,,
\ee
where $\beta\in\{1,\alpha\}$. It is enough to consider only the case $\beta=1$.
Clearly $I(x)$ is bounded for $|x|\le 1$. To estimate $I(x)$ when $|x|$ is large,
let us write $x=te$ with $|e|=1$ and make the substitution $y=tz$ in~\eqref{integral} (with $\beta=1$). 
We obtain
\be
\label{integral2}
I(te)=t^{-\aa}\int_{\rt}{\frac{dz}{|e-z|^2(t^{-1}+|z|)^{\aa+1}}}\le t^{-\aa}\int_{\rt}{\frac{dz}{|e-z|^2|z|^{\aa+1}}}\,\,.
\ee
Since we assume $\aa\in(1,2)$, the last integral is bounded, and we see that 
\be
\label{integral3}
t^\aa I(te)\le c(\aa)\,\,.
\ee
Combining this estimate with the estimate of $I(x)$ for $|x|\le 1$ we see that
\be
\label{integral4}
(1+|x|)^\aa I(x)\le c(\aa)\qquad \mbox{for all $x\in\rt\,\,.$}
\ee
This completes the proof of estimates~\eqref{estimate1} and~\eqref{estimate2}.
\end{myproof}

We have shown that the linear operator
$$
{T_{\tuu}} \colon v\to G*\div(\tuu \otimes v + v\otimes\tuu)
$$ 
is continuous from $\xa$ to $\xa$, and its norm is bounded by $c(\aa,M)|b|$. Also, we have shown that the bi-linear
operator
$$
B\colon (v,w)\to G*\div(v\otimes w)
$$
is continuous from $\xa\times\xa\to\xa$, with the bound
$$
\norm{B(v,w)}\le c(\aa)\norm{v}\norm{w}\,\,.
$$
We let $V=G*(f-\tf)$ and re-write equation~\eqref{maineq} as
\be
\label{opereq}
v+T_{\tuu}(v)+B(v,v)=V\,\,.
\ee
Since $\int_{\rt}(f-\tf)=0$, we have $V=O(|x|^{-2})$ as $|x|\to\infty$. 
Standard perturbation arguments (such as the Implicit Function Theorem) now imply that equation~\eqref{opereq}
has a solution $v$ when $V$ is sufficiently small in $\xa$. (A simple sufficient condition for that is
that, in addition to  $\int_{\rt}(f-\tf)=0$ and the restriction on the support on $f-\tf$, 
the field $f-\tilde F$ be small in $L^{\frac32 + \delta}$ with some $\delta>0$.) 
 Moreover, the solution is unique in some small
ball in $\xa$ (centered at the origin). These statements can be made more quantitative if we use the
special form of the perturbation (namely that it is quadratic in~$v$). For example, one can use the following
folklore lemma:
\begin{Lemma}
\label{abstractlemma}
Let $X$ be a Banach space. Let $T\colon X\to X$ be linear with $||Tx||\le \ve ||x||$ for all $x\in X$, 
and let $B\colon X\times X\to X$ be bilinear with $||B(x_1,x_2)||\le c||x_1||\,||x_2||$
for all $x_1,x_2\in X$. Let $y\in X$ with $||y||<{\frac{(1-\ve)^2}{4c}}$. Let $0<\xi_1<\xi_2$ be the
two roots of the equation $\xi=||y||+\ve\xi+c\xi^2$, i.\ e.\ $\xi_{1,2}={\frac{(1-\ve)\mp\sqrt{(1-\ve)^2-4c||y||}}{2c}}$. Then the equation
\be
\label{abstracteq}
x+Tx+B(x,x)=y
\ee
has a solution $\bar x$ satisfying $||\bar x||\le\xi_1$. Moreover, the solution $\bar x$ is unique in the open ball $\{x\in X, ||x||<\xi_2\}$. 
\end{Lemma}

\begin{myproof}
The proof is standard and we include it for the convenience of the reader. Consider the map $F(x)=y-Tx-B(x,x)$.
We have $||F(x)||\le ||y||+\ve||x||+c||x||^2$ which shows that for $\xi_1<||x||<\xi_2$ we have
$||F(x)||<||x||$ and that, for any $\delta>0$,  the iterates $F(x), F^2(x)=F(F(x)),\dots, F^k(x),\dots$ enter the ball
of radius $\xi_1+\delta$ after finitely many steps. At the same time, we have
$||F(x_1)-F(x_2)||\le \ve||x_1-x_2||+c||x_1-x_2||(||x_1||+||x_2||)$ which shows that $F$ is a contraction of any 
closed ball of radius $\xi\in[\xi_1,{\frac{\xi_1+\xi_2}2})$. 
\end{myproof}

\begin{myproofb} 
Let us first assume that the ``no outflow to infinity" condition~\eqref{nooutflow1} is satisfied.
In this case the statement of Theorem~\ref{mainth} is a direct consequence of the construction of the extensions in Section~\ref{preliminaries} and
Lemmata~\ref{elementaryest} and~\ref{abstractlemma}. Note that we not only need existence and uniqueness for $v$ in~\eqref{maineq}, but we also
need uniqueness (with smallness assumptions) for $u$ in~\eqref{nseagain}. The uniqueness of $u$ in our situation is well-known (see e.\ g.\
\cite{Finn} or \cite{Galdi}), and can also be easily proved from Lemma~\ref{abstractlemma} and (an obvious modification
of) Lemma~\ref{elementaryest}.

The situation when we have some outflow to infinity can be handled by a standard method of using the canonical
outflow field $\frac{x\,\,}{|x|^3}$, see for example~\cite{Finn}, Section 2.2, \cite{Galdi}, Chapter IX, or \cite{NazarovPileckas}, Remark 3.2. 
Assume $R_0=1$ without loss of generality. Let $a$ be the multiple of the vector field $\frac{x\,\,}{|x|^3}$ which has the same 
outflow as $u$. Note that $a$ satisfies
the Navier-Stokes equation~\eqref{nse} in $\rt\setminus\{0\}$ with the associated pressure field $\pi_a=-{\frac12}|a|^2$.
Let us write $u=a+w$ and $p=\pi_a+p_w.$ The field $w$ satisfies the no outflow condition, and we can extend it to a div-free field $\tilde w$
with the control similar to~\eqref{extestu}. We can also regularize $a$ and $\pi_a$ in $B_1$ (while not changing them outside $B_1$)
so that estimates similar to~\eqref{extestu} and~\eqref{extestp} are satisfied. Let us denote be $\tilde a$ and $\tilde \pi_a$
these regularized functions. Finally, we extend $p_w$ to $\tilde p_w$ with control similar to~\eqref{extestp}.
Let $\tilde u=\tilde a +\tilde w$ and $\tilde p=\tilde \pi_a +\tilde p_w$. (Note that $\tilde u$ is not div-free in $B_1$.)
Let $\tilde f=\div \tilde T$, where $\tilde T=\tilde T(\tilde u,\tilde p)$ is given by~\eqref{mflux} with $u,p$ replaced by $\tilde u, \tilde p$.
We note that the vector $b$ given by~\eqref{formulab} can also be expressed as  $b=\int_{\rt}\tilde f\,.$
We will now search a div-free vector field $z$ and a function $p_z$ satisfying $\div \tilde T(\tilde a + z,\tilde\pi_a+p_z)=\tilde f$.
We seek $z,p_z$ in the form $z={\tilde U}^b+v,\quad p_z={\tilde P}^b+q$, where ${\tilde U}^b, {\tilde P^b}$ are the regularizations
of the Landau solutions $U^b, P^b$ constructed in Section~\ref{landau}. It is now easy to check that the perturbation
theory of Section~\ref{perturbation} gives the required solution.

\end{myproofb}

\begin{Remark}
\label{remark1}
{\rm  The borderline space $\xa$ in which a more sophisticated perturbation analysis might possibly work is the 
space $X_2$. (This corresponds to the naturally expected decay $O(|x|^{-2})$ for $v$.) However, a perturbation analysis 
in $X_2$ cannot be based only on the decay properties of $\tuu$ (as was the case with our simpler analysis for
$\aa<2$). To see this, let $\ve\in(0,1)$ and consider the equation 
\be
\label{example}
-\Delta u + {\frac{\ve(n-\ve)}{1-\ve}}\div({\frac{x}{\,\,|x|^2}}\,\,u)=f(x)\qquad\mbox{in $\rr^n\,.$}
\ee
One can check by direct calculation that when $f=0$ the function $x_1|x|^{-n+\ve}$ is a solution of this equation 
away from the the origin. Let $\eta$ be a smooth function in $\rr^n$ which vanishes in the unit ball
and is equal to $1$ outside of the ball of radius $2$. An easy calculation shows that
the function $u=\eta \,x_1|x|^{-n+\ve}$ satisfies~\eqref{example} with $\int_{\rt} f = 0\,\,$. Moreover, one can change the coefficients
of the equation in the unit ball so that they become smooth.

To get results in the space $X_2$, one would probably have to prove optimal decay estimates for the linear equation
\be
\label{lmaineq}
\left.\begin{array}{rcl}
-\Delta v +\tuu\nabla v +v\nabla \tuu +\nabla q & = & f-\tf\,\, , \\
\div v & = & 0
\end{array}\right\}\qquad \mbox{in $\rt$}
\ee
by a non-perturbative approach, and then treat the quadratic term in~\eqref{maineq} perturbatively.
The above example shows that to get the optimal decay $O(|x|^{-2})$, one would need to use more
information about $\tuu$ than just its decay properties at $\infty$.
We conjecture that for large $|x|$ the perturbation $v$ from the Landau solution $\ub$ indeed
has the decay $v(x)=O(|x|^{-2})$, at least for small data.

}
\end{Remark}

\end{document}